\documentclass [a4paper] {article}

\usepackage[utf8]{inputenc}

\usepackage[french]{babel}

\usepackage{amssymb,amsmath, amsfonts}

\usepackage{amsthm}

\usepackage{hyperref}

\usepackage{MnSymbol}

\newtheorem{prop} {Proposition} 
\newtheorem{lm} [prop]{Lemme} 
\newtheorem{thm} [prop] {Théorème}

\theoremstyle{definition}
\newtheorem{df}{Définition} 
\newtheorem*{df*}{Définition}

\theoremstyle{remark}
\newtheorem{rmq}{Remarque}

\newcommand*{\twosquig}{%
\mathrel{\vcenter{\offinterlineskip
\hbox{$\rightsquigarrow$}\vskip-.08ex\hbox{$\rightsquigarrow$}}}}

\author{Stéphane \textsc{Dugowson}
\footnote{s.dugowson@gmail.com}
}

\title {{\Large {Dynamiques sous-catégoriques ouvertes en interaction}\\(définitions et théorème de stabilité)}}
\date{9 août 2015} 

\begin{document}

\maketitle

\noindent \textbf{Abstract}. \textsf{Open sub-categorical dynamics in interaction (definitions and stability theorem).} The aim of this paper is to define what we  call \emph{open sub-categorical dynamics}, their interactions and the sub-categorical dynamics produced by those interactions, thanks to the stability theorem we prove here and which motivates all this study.

\mbox{}

\noindent \textit{Keywords}:  Interaction. Open dynamics. Sub-categorical dynamics. Relations. Multiple binary relations. Connectivity structures. Determinism. Indeterminism.

\mbox{}

\noindent\textbf{Résumé}. En définissant les dynamiques sous-catégoriques ouvertes et leurs interactions, cet article se situe dans le prolongement immédiat de notre article précédent consacré aux dynamiques graphiques ouvertes. Les dynamiques sous-catégoriques constituent une généralisation des dynamiques catégoriques, présentant sur ces dernières l'avantage d'une certaine stabilité, objet du \og théo\-rème de stabilité\fg\, que nous prouvons à la fin de l'article.

\mbox{}

\noindent \textit{Mots clés}:  Interaction.  Dynamiques ouvertes. Dynamiques sous-catégoriques. Relations. Relations binaires multiples. Structures connectives. Déterminisme. Indéterminisme.

\mbox{}

\noindent Mathematics Subject Classification 2010: 18A10. 37B55. 54H20.

\section{Introduction}

Les \emph{dynamiques sous-catégoriques} généralisent les dynamiques catégoriques définies dans \cite{Dugowson:201112} et \cite{Dugowson:201203}, et les \emph{dynamiques sous-catégoriques ouvertes} généralisent les dynamiques catégoriques ouvertes présentées dans notre conférence sur les \og Influences dynamiques catégoriques\fg\, au cours du colloque \og Topologie, effectivité, interactivité\fg\, organisé à Supméca le 21 mai 2013 \cite{Dugowson:20130521}. Cette généralisation est, en effet, rendue nécessaire par le fait que  la catégoricité des dynamiques catégoriques ouvertes se révèle instable dans leurs interactions\footnote{Nous y reviendrons dans notre conclusion. Voir en particulier la note \ref{footnote principe exemple cat instable}.}. Pour obtenir la stabilité que nous recherchions --- et qu'exprime ici le théorème \ref{thm stabilite sc} donné page \pageref{thm stabilite sc} --- nous avons dû dans un premier temps nous replier sur la notion assez restreinte mais déjà complexe de \emph{dynamique graphique}, ce que nous avons fait avec notre texte \og Interaction des dynamiques graphiques ouvertes\fg\, \cite{Dugowson:20150807} qui constitue le socle sur lequel toutes les définitions qui suivent seront données. \`{A} quelques nuances près
, nous reprenons ici toutes les notations et définitions introduites dans \cite{Dugowson:20150807}, ainsi que celles  de  \cite{Dugowson:201112}  et \cite{Dugowson:201203}. Ainsi,
\begin{itemize}
\item dans tout l'article, $\mathbf{C}$ désigne une petite catégorie,
\item pour tout objet $A$ d'une catégorie, $Id_A$ désigne le morphisme identité,
\item $Gr(\mathbf{C})$ désigne le graphe associé à $\mathbf{C}$ en \og oubliant\fg\, la composition des flèches,
\item $\overrightarrow{\mathbf{C}}$ désigne l'ensemble des flèches de $\mathbf{C}$ et $\dot{\mathbf{C}}$ l'ensemble de ses objets (\emph{i.e.} les sommets du graphe $Gr(\mathbf{C})$),
\item pour toute flèche $e:A\rightarrow B$ d'une catégorie (ou pour toute arête d'un graphe), $dom(e)$ désigne sa source (ou domaine) $A$ et $cod(e)$ désigne son but (ou codomaine) $B$,
\item $\mathbf{P}$ désigne la catégorie dont les objets sont les ensembles et dont les flèches sont les relations binaires entre ensembles, ces relations étant vues en tant que \og transitions non déterministes\fg, dont la composition est notée $\odot$,
\item Pour tout ensemble $M$, $\mathbf{P}^{\underrightarrow{{M}}}$ désigne la catégorie dont les objets sont les ensembles et dont les flèches sont les familles indexées par $M$ de transitions non déterministes,
\item  si $S$ et $T$ sont deux ensembles et  $u:S\leadsto T$ et $v:S\rightsquigarrow T$ deux transitions, nous écrirons 
\[u\subset v\] pour exprimer le fait que 
\[\forall a\in S, u(a)\subset v(a),\]
\item etc...
\end{itemize}

Par ailleurs, pour toute transition $u:S\leadsto T$ et tout élément $a\in S$, si $u(a)$ est un singleton : $u(a)=\{b\}$ et s'il n'y a pas d'ambiguïté, on commettra souvent l'abus d'écriture consistant à écrire $u(a)=b$. En particulier, pour une fonction $f:S\rightarrow T$ non partout définie, nous pourrons noter $f(a)=b$ si $f$ est définie en $a$ et y prend la valeur $b$, et $f(a)=\emptyset$ si $f$ n'est pas définie en $a$.

\section{Dynamiques sous-catégoriques ouvertes}
\label{section dynamiques souscat ouvertes}

\subsection{Dynamiques sous-catégoriques}

\subsubsection{Définition des dynamiques sous-catégoriques}

\begin{df}\label{df dyscat} Une \emph{dynamique sous-catégorique} $\alpha$ consiste en la donnée d'une petite catégorie $\mathbf{C}$, appelée \emph{moteur} de la dynamique $\alpha$, et d'une dynamique graphique\footnote{Voir \cite{Dugowson:20150807}, définition 6.} $Gr(\alpha):Gr(\mathbf{C})\rightarrow Gr(\mathbf{P})$ telle que 
\begin{itemize}
\item pour tout objet $A\in\dot{\mathbf{C}}$, on a \[(Id_A)^{Gr(\alpha)}\subset Id_{A^{Gr(\alpha)}},\]
\item pour tout couple $(e,f)\in \overrightarrow{\mathbf{C}}^2$ de flèches composables --- \emph{i.e.} telles que $cod(e)=dom(f)$ --- on a
\[
(f\circ e)^{Gr(\alpha)} \subset f^{Gr(\alpha)}\odot e^{Gr(\alpha)}.
\]
\end{itemize}
\end{df}

On écrira $\alpha:\mathbf{C}\rightharpoondown \mathbf{P}$ pour exprimer que $\alpha$ est une dynamique sous-catégorique de moteur $\mathbf{C}$. Par ailleurs, la dynamique graphique ${Gr(\alpha)}$ sera le plus souvent simplement désignée  par la même lettre $\alpha$, de sorte que les relations de la définition \ref{df dyscat} ci-dessus peuvent s'écrire
\[(Id_A)^{\alpha}\subset Id_{A^{\alpha}},\] et
\[
(f\circ e)^{\alpha} \subset f^{\alpha}\odot e^{\alpha}.
\]

Les flèches de la catégorie $\mathbf{C}$ seront également appelées les \emph{$\mathbf{C}$-écoulements}, les \emph{écoulements de $\mathbf{C}$}, ou simplement, s'il n'y a pas d'ambiguïté, les \emph{écoulements}.

L'ensemble $st(\alpha)$ des états d'une dynamique sous-catégorique $\alpha$, le type $typ(s)$ d'un état $s\in st(\alpha)$, les notions de déterminisme et de quasi-déterminisme appliquées à $\alpha$ sont respectivement définis par les notions analogues\footnote{Voir la section §\textbf{2.1.1} de \cite{Dugowson:20150807}.} appliquées à la dynamique graphique $Gr(\alpha)$.

\subsubsection{Catégories de  dynamiques sous-catégoriques}

Les définitions, données en section §\textbf{2.1.2} de \cite{Dugowson:20150807}, des dynamorphismes entre dynamiques graphiques, qu'on se restreigne au cas de dynamiques ayant même moteur ou qu'on considère le cas général, s'étendent aux dynamiques sous-catégoriques. On obtient ainsi les définitions suivantes des catégories $\mathbf{DySC_{(C)}}$ et $\mathbf{DySC}$ :

\begin{df}\label{df dynamorphismes dyscat C} Étant données deux dynamiques sous-catégoriques $\alpha : \mathbf{C}\rightharpoondown  \mathbf{P}$ et $\beta : \mathbf{C}\rightharpoondown  \mathbf{P}$ de \emph{même} moteur $\mathbf{C}$, un \emph{dynamorphisme $\delta:\alpha\looparrowright_{\mathsf{SC}}\beta$ sous-catégorique sur $\mathbf{C}$} est un dynamorphisme   $\delta:Gr(\alpha)\looparrowright Gr(\beta)$ entre dynamiques graphiques de même moteur $Gr(\mathbf{C})$. On notera $\mathbf{DySC_{(C)}}$ la catégorie dont les objets sont les dynamiques sous-catégoriques de moteur $\mathbf{C}$ et dont les flèches sont les dynamorphismes sous-catégoriques sur $\mathbf{C}$.
\end{df}

Autrement dit, dans $\mathbf{DySC_{(C)}}$, un morphisme $\delta:\alpha\looparrowright_{\mathsf{SC}}\beta$ est la donnée d'une transition $\delta:st(\alpha)\rightsquigarrow st(\beta)$ telle que
\[\forall A\in\dot{\mathbf{C}},\forall a\in A^\alpha, \delta(a)\subset A^\beta,\] et
\[\forall (f:A\rightarrow B)\in\overrightarrow{\mathbf{C}},\forall  a\in A^\alpha, (\delta\odot f^\alpha)(a)\subset (f^\beta\odot\delta)(a).\]
 
La petite catégorie $\mathbf{C}$ étant donnée, la catégorie $\mathbf{DySC_{(C)}}$ s'identifie à une sous-catégorie pleine de la catégorie des dynamiques graphiques de moteur le graphe $Gr(\mathbf{C})$. Par conséquent, il n'y aura aucune ambiguïté à écrire ${\delta:\alpha\looparrowright\beta}$ au lieu de $\delta:\alpha\looparrowright_{\mathsf{SC}}\beta$ pour exprimer que $\delta$ est un dynamorphisme sous-catégorique entre dynamiques sous-catégoriques de même moteur.

\begin{df}\label{df dynamorphismes dyscat} Étant données deux dynamiques sous-catégoriques $\alpha : \mathbf{C}\rightharpoondown  \mathbf{P}$ et $\beta : \mathbf{D}\rightharpoondown  \mathbf{P}$,
 un \emph{dynamorphisme sous-catégorique} $(\Delta,\delta):\alpha\looparrowright_{\mathsf{SC}}\beta$
  est un dynamorphisme de dynamiques graphiques\footnote{Voir la section §\textbf{2.1.2} de \cite{Dugowson:20150807}.} $(\Delta,\delta):Gr(\alpha)\looparrowright Gr(\beta)$ 
  tel que $\Delta:\mathbf{C}\rightarrow\mathbf{D}$ soit un foncteur.  
  On notera $\mathbf{DySC}$ la catégorie dont les objets sont les dynamiques sous-catégoriques   et dont les flèches sont les dynamorphismes sous-catégoriques.
\end{df}

S'il n'y a pas d'ambiguïté, on dira que $(\Delta,\delta)$ est un dynamorphisme de $\alpha$ vers $\beta$, et on écrira  ${(\Delta,\delta):\alpha\looparrowright\beta}$ (au lieu de $(\Delta,\delta):\alpha\looparrowright_{\mathsf{SC}}\beta$)  pour exprimer que $(\Delta,\delta)$ est un dynamorphisme sous-catégorique entre les dynamiques sous-catégoriques $\alpha$ et $\beta$.

\subsubsection{Foncteurs d'oubli $Gr:\mathbf{DySC_{(C)}}\rightarrow \mathbf{DyG}_{Gr(\mathbf{C})}$}

Notons $\mathbf{DyG}_{(\mathbf{G})}$ (respectivement $\mathbf{DyG}$) la catégorie des dynamiques graphiques de moteur le graphe $\mathbf{G}$ (respectivement la catégorie de toutes les dynamiques graphiques). On notera $Gr$ le foncteur d'oubli $\mathbf{DySC}\rightarrow \mathbf{DyG}$ qui à toute dynamique sous-catégorique $\alpha$ associe la dynamique graphique $Gr(\alpha)$ et à tout dynamorphisme sous-catégorique $\alpha \looparrowright \beta$ associe lui-même en tant que dynamorphisme graphique $Gr(\alpha)\looparrowright Gr(\beta)$. Par restriction aux dynamiques sous-catégoriques de moteur $\mathbf{C}$ et aux dynamorphismes sous-catégoriques sur $\mathbf{C}$, on obtient un foncteur d'oubli, encore noté $Gr$, de $\mathbf{DySC_{(C)}}$ vers $\mathbf{DyG}_{Gr(\mathbf{C})}$. Pour souligner l'application de $Gr$, on pourra si nécessaire la faire porter sur chacun des termes de l'expression concernée. Ainsi, à un dynamorphisme sous-catégorique
\[
(\Delta,\delta):
(\alpha:\mathbf{C}\rightarrow \mathbf{P})
\looparrowright
(\beta:\mathbf{D}\rightarrow \mathbf{P}),\] le foncteur d'oubli associe le dynamorphisme graphique
\[
(Gr(\Delta),Gr(\delta)):
(Gr(\alpha):Gr(\mathbf{C})\rightarrow Gr(\mathbf{P}))
\looparrowright
(Gr(\beta):Gr(\mathbf{D})\rightarrow Gr(\mathbf{P})).\] 

\subsubsection{Dynamiques sous-catégoriques propres}\label{subsubsec Dyna SCat propres}

Soit $\alpha:\mathbf{C}\rightharpoondown \mathbf{P}$ une dynamique sous-catégorique. Pour tout état $a\in A^\alpha$ de type $A\in\dot{\mathbf{C}}$, on a soit $(Id_A)^\alpha(a)=a$, soit $(Id_A)^\alpha(a)=\emptyset$. Dans le  cas où $(Id_A)^\alpha(a)=\emptyset$, nous dirons que, pour $\alpha$, $a$ est \emph{en-dehors du coup}.

Pour un tel état $a$ en-dehors du coup, on a pour toute flèche $(f:A\rightarrow B)\in\overrightarrow{\mathbf{C}}$
\[f^\alpha(a)=(f\circ Id_A)^\alpha(a)\subset f^\alpha \odot (Id_A)^\alpha(a)=f(\emptyset)=\emptyset,\] et de même, pour toute flèche $g:B\rightarrow A$ et tout $b\in B^\alpha$
\[a\in g^\alpha(b)\Rightarrow a\in (Id_A \circ g)^\alpha (b) \subset (Id_A)^\alpha \odot g^\alpha (b)= \bigcup_{a'\in g^\alpha (b)} (Id_A)^\alpha(a'), \]
Mais 
\[(Id_A)^\alpha(a)=\emptyset \Rightarrow \bigcup_{a'\in g^\alpha (b)} (Id_A)^\alpha(a')\subset g^\alpha (b)\setminus \{a\} 
,\] de sorte que $a$ appartient à un ensemble auquel il n'appartient pas : on en déduit par l'absurde qu'en fait 
\[
\forall (g:B\rightarrow A), \forall b\in B^\alpha, a\notin g^\alpha(b). 
\]

Ainsi, un état en-dehors du coup a pour toute transition de la dynamique considérée une image vide, et n'est dans l'image d'aucune transition, pour aucun état.

On vérifie facilement que l'on peut alors retirer fonctoriellement tous les états en dehors du coup de toute dynamique sous-catégorique $\alpha$ pour obtenir une dynamique $\check{\alpha}$ dont tous les états soient dans le coup, 
de sorte qu'à tout dynamorphisme $\delta:\alpha\looparrowright\beta$ entre de telles dynamiques sous-catégoriques se trouve associé un dynamorphisme $\check{\delta}:\check{\alpha}\looparrowright\check{\beta}$. Nous appellerons dynamiques sous-catégoriques propres les dynamiques sous-catégoriques ainsi \og nettoyées\fg\,:
\begin{df}\label{df DySCP} Une \emph{dynamique sous-catégorique propre} $\alpha$ est une dynamique sous-catégorique dont tous les états sont dans le coup, autrement dit telle que
\[\forall A\in\dot{\mathbf{C}}, (Id_{A})^\alpha =Id_{A^\alpha}.\]
\end{df}

On définit la catégorie des dynamiques sous-catégoriques propres comme la sous-catégorie pleine de la catégorie des dynamiques sous-catégoriques obtenue en se restreignant aux dynamiques sous-catégoriques propres.

\subsubsection{Dynamiques catégoriques}

\begin{df}
Une \emph{dynamique catégorique}  de moteur la petite catégorie $\mathbf{C}$ est une dynamique sous-catégorique propre $\alpha$ vérifiant l'égalité
\[
(f\circ e)^{\alpha} = f^{\alpha}\odot e^{\alpha}.
\]
\end{df}

\begin{rmq}
La définition ci-dessus équivaut à celle d'une dynamique catégorique \emph{propre} donnée dans \cite{Dugowson:201112} et \cite{Dugowson:201203}. On fera toutefois attention au fait que le qualificatif \emph{propre} n'y a plus la même signification : dans les textes cités, il se rapporte aux dynamiques catégoriques $\alpha$ pour lesquelles $S\neq T$ $ \Rightarrow$ $ S^\alpha\cap T^\alpha= \emptyset$, condition qui est automatiquement satisfaite ici (car ayant été incorporée dans la définition des dynamiques graphiques), tandis qu'il se rapporte maintenant aux dynamiques sous-catégoriques  qui vérifient l'égalité $(Id_A)^\alpha=Id_{A^\alpha}$.
\end{rmq}

\begin{prop}\label{prop dysc deterministe implique categorique}
Une dynamique sous-catégorique déterministe est nécessairement catégorique.
\end{prop}
\paragraph{Preuve.} Si $\alpha$ est déterministe, les deux membres des inclusions de la forme
\[(g\circ f)^{\alpha}(a) \subset (g^{\alpha}\odot f^{\alpha})(a)\] sont des singletons, il y a donc égalité.
\begin{flushright}$\square$\end{flushright} 

\subsubsection{Horloge, successions, réalisations}
\label{subsubsec horloge cat success real}

\begin{df}\label{df horloge cat} Une \emph{horloge} de moteur la petite catégorie $\mathbf{C}$ est une dynamique sous-catégorique déterministe sur $\mathbf{C}$.
\end{df}

D'après la proposition \ref{prop dysc deterministe implique categorique}, une horloge est une dynamique catégorique, de sorte que la définition \ref{df horloge cat} ci-dessus est équivalente à celle donnée dans \cite{Dugowson:201112} et \cite{Dugowson:201203}.  La partie graphique $Gr(h)$ d'une horloge $h$ est une horloge sur le graphe $Gr(\mathbf{C})$, et conformément à la définition d'une horloge graphique, les états de $h$ seront également appelés des \emph{instants}.

Rappelons également que la relation de succession entre instants d'une telle horloge catégorique $h$, 
\[s\leq t \Leftrightarrow \exists e\in \overrightarrow{\mathbf{C}}, e^h(s)=t.\] relation qui n'avait aucune propriété particulière dans le cas graphique, est à présent une relation de pré-ordre.

Prolongeant la définition 9 de \cite{Dugowson:20150807}, on pose :

\begin{df}\label{df realisation de DySC pour horloge} Étant donnée $\mathbf{C}$ une petite catégorie, $\alpha$ une dynamique sous-catégorique de moteur $\mathbf{C}$ et $h$ une horloge de même moteur, une \emph{$h$-réalisation} (ou \emph{$h$-solution}) $\mathfrak{s}:h\looparrowright \alpha$  de $\alpha$ est une $Gr(h)$-réalisation de $Gr(\alpha)$.
\end{df}

\begin{rmq}\label{rmq meme real pour dysc et graph} Soulignons que les réalisations d'une dynamique sous-catégoriques sont les mêmes que celles de cette même dynamique considérée en tant que dynamique graphique. 
\end{rmq}

Autrement dit, une $h$-réalisation de la dynamique sous-catégorique $\alpha:\mathbf{C}\rightharpoondown\mathbf{P}$ est une fonction $\mathfrak{s}:st(h)\dashrightarrow st(\alpha)$ définie sur une partie de l'ensemble des instants de $h$, vérifiant 
\[
\forall (f:S\rightarrow T)\in\overrightarrow{\mathbf{C}}, 
\forall s\in S^h,
\mathfrak{s}(f^h(s))\subset f^\alpha(\mathfrak{s}(s)),
\]
ce qui signifie que 
\begin{itemize}
\item soit $\mathfrak{s}$ n'est pas définie à l'instant $s$, auquel cas elle ne l'est pas non plus à l'instant ultérieur $f^h(s)$ puisque
dans ce cas
\[\mathfrak{s}(f^h(s))\subset f^\alpha(\mathfrak{s}(s))=\emptyset
\Rightarrow \mathfrak{s}(f^h(s))=\emptyset
,\]
\item soit $\mathfrak{s}$ est défini à l'instant $s$  et  n'est pas définie à l'instant ultérieur $f^h(s)$, et dans ce cas on a 
\[\emptyset=\mathfrak{s}(f^h(s))\subsetneqq f^\alpha(\mathfrak{s}(s)),\]
\item soit $\mathfrak{s}$ est  définie à l'instant  $f^h(s)$, et dans ce cas elle est nécessairement définie à l'instant $s$ avec\footnote{\label{foot ecriture transition quasidet}
Rappelons la convention que nous avons adoptée d'écrire $f(a)=b$ lorsque $f$ est une transition vérifiant $f(a)=\{b\}$ pour certains états $a$ et $b$.}
\[\mathfrak{s}(f^h(s))\in f^\alpha(\mathfrak{s}(s)).\]
\end{itemize}

Bien entendu, dans le cas où $\alpha$ est une dynamique \emph{déterministe}, le troisième cas ci-dessus s'écrit \[\mathfrak{s}(f^h(s))= f^\alpha(\mathfrak{s}(s)).\]

\subsubsection{Union de dynamiques sous-catégoriques sur $\mathbf{C}$}

\paragraph{Cas des dynamiques graphiques de moteur $\mathbf{G}$.} Étant donnée $(\alpha_\nu:\mathbf{G}\rightarrow Gr(\mathbf{P}))_{\nu\in N}$ une famille indexée par un ensemble quelconque $N$ de dynamiques graphiques de même moteur un graphe $\mathbf{G}$, on appelle union des dynamiques de cette famille, et on note $\bigcup_{\nu\in N} \alpha_\nu$, la dynamique graphique $\alpha$ de moteur $\mathbf{G}$ telle que
\begin{itemize}
\item pour tout sommet $S\in \dot{\mathbf{G}}$, 
\[S^\alpha=\bigcup_{\nu\in N} S^{\alpha_\nu},\]
\item pour toute arête $(f:S\rightarrow T)\in\overrightarrow{\mathbf{C}}$ et tout état $s\in S^\alpha$, 
\[f^\alpha(s)=\bigcup_{\nu\in N} f^{\alpha_\nu}(s),\]où l'on convient de prolonger $f^{\alpha_\nu}$ à $S^\alpha$ en posant $f^{\alpha_\nu}(s)=\emptyset$ si $s\notin S^{\alpha_\nu}$.
\end{itemize}

\begin{rmq}
L'union de la famille vide de dynamiques graphiques sur $\mathbf{G}$ est la dynamique vide sur $\mathbf{G}$ : l'ensemble des états est vide.
\end{rmq}

\begin{prop}\label{prop union dyna sous-cat} Étant donnée $(\alpha_\nu:\mathbf{C}\rightharpoondown \mathbf{P})_{\nu\in N}$ une famille indexée par un ensemble $N$ de dynamiques sous-catégoriques de même moteur une petite catégorie $\mathbf{C}$, alors $\bigcup_{\nu\in N} Gr(\alpha_\nu)$,  la dynamique graphique définie par l'union des dynamiques $\alpha_\nu$ vues comme dynamiques graphiques, est elle-même une dynamique sous-catégorique. Si, en outre, toutes les dynamiques sous-catégoriques $\alpha_\nu$ sont propres, il en va de même de leur union.
\end{prop}
\paragraph{Preuve.} Posons $\alpha=\bigcup_{\nu\in N} Gr(\alpha_\nu)$. Soit $S\in\dot{\mathbf{C}}$. On a
\[(Id_S)^\alpha=\bigcup_{\nu\in N} {(Id_S)}^{\alpha_\nu}\subset Id_{S^\alpha}.\] Dans le cas où tous les $\alpha_\nu$ sont propres, il est immédiat qu'on a en fait $(Id_S)^\alpha=Id_{S^\alpha}$.
Soit maintenant $(S\stackrel{f}{\rightarrow}T)$ et $(T\stackrel{g}{\rightarrow}U)$ deux flèches composables de $\mathbf{C}$. On a
\[
(g\circ f)^\alpha
=
\bigcup_{\nu\in N}(g\circ f)^{\alpha_\nu}
\subset
\bigcup_{\nu\in N}(g^{\alpha_\nu}\odot f^{\alpha_\nu})
\subset
\bigcup_{\nu\in N}(g^{\alpha}\odot f^{\alpha})
=
(g^{\alpha}\odot f^{\alpha}).
\] 
\begin{flushright}$\square$\end{flushright} 

Munissant la classe des dynamiques sous-catégoriques sur $\mathbf{C}$ de la relation d'ordre $\subset$ définie par 
\[\alpha\subset\beta \Leftrightarrow 
\left\lbrace
\begin{tabular}{l}
$\forall S\in\dot{G}, S^\alpha\subset S^\beta$, \\ 
$\forall (S\stackrel{f}{\rightarrow}T)\in\overrightarrow{\mathbf{C}}, f^\alpha\subset f^\beta$, \\ 
\end{tabular} 
\right.
\]
on déduit de la proposition \ref{prop union dyna sous-cat} que, pour toute dynamique sous-catégorique $\alpha:\mathbf{C}\rightharpoonup \mathbf{P}$, l'ensemble des dynamiques sous-catégoriques de même moteur qui sont incluses dans $\alpha$ constitue  un treillis complet avec l'union pour $\sup$. En particulier, à toute dynamique graphique $\gamma$ de moteur $Gr(\mathbf{C})$, on peut associer la plus grande dynamique sous-catégorique $\stackrel{\circ}{\gamma}$ qui soit plus petite que $\gamma$ :
\[
\stackrel{\circ}{\gamma}=\bigcup_{\beta:\mathbf{C}\rightharpoonup\mathbf{P},\\ Gr(\beta)\subset \gamma}\beta.
\]

Un autre corolaire immédiat de la proposition \ref{prop union dyna sous-cat} est le fait que l'union d'une famille (indexée par un ensemble) de dynamiques catégoriques sur $\mathbf{C}$ est une dynamique sous-catégorique propre sur $\mathbf{C}$. Par contre, en général, l'union de dynamiques catégoriques n'est pas catégorique en général. Par exemple, pour $\mathbf{C}$ engendrée par les flèches $1\rightarrow 2\rightarrow 3$, les deux dynamiques catégoriques quasi-déterministes $\alpha_1$ et $\alpha_2$ ayant mêmes ensembles d'états
\[1^\alpha=\{a_1\},\quad 2^\alpha=\{a_2,a'_2\},\quad 3^\alpha=\{a_3,a'_3\} \]
et dont les transitions sont engendrées respectivement par
\[ (1\rightarrow 2)^{\alpha_1}(a_1)=a_2\quad \mathrm{et}\quad 
(2\rightarrow 3)^{\alpha_1}(a_2)=a_3
\]
et
\[ (1\rightarrow 2)^{\alpha_2}(a_1)=a'_2,\quad
(2\rightarrow 3)^{\alpha_2}(a'_2)=a_3
\quad \mathrm{et}\quad 
(2\rightarrow 3)^{\alpha_2}(a_2)=a'_3
\]
ont pour union $\alpha=\alpha_1\cup \alpha_2$ vérifiant
\[
(1\rightarrow 3)^\alpha(a_1)
=
\{a_3\}
\subsetneqq
\{a_3,a'_3\}=
\left((2\rightarrow 3)^{\alpha}\odot(1\rightarrow 2)^{\alpha}\right) (a_1).
\]

\begin{rmq}
Symétriquement, l'intersection de dynamiques catégoriques n'est pas, en général, une dynamique catégorique, mais seulement une dynamique \emph{extra-catégorique} propre, c'est-à-dire vérifiant les inclusions de la forme
\[(f\circ e)^{\alpha} \supset f^{\alpha}\odot e^{\alpha}.\] Un exemple très simple de cette situation est considéré dans notre conférence \cite{Dugowson:20150506}.
\end{rmq}

\subsubsection{Récapitulatif pour les dynamiques sous-catégoriques}

En résumé, les dynamiques sous-catégoriques et les notions qui leur sont relatives sont définies par le fait que les foncteurs d'oubli $Gr$ les envoient sur les dynamiques graphiques et les notions associées, avec les contraintes supplémentaires suivantes : une dynamique sous-catégorique doit vérifier
\[(f\circ g)^\alpha \subset f^\alpha \odot g^\alpha,\]
et
\[(Id_A)^\alpha \subset Id_{A^\alpha},\]
une dynamique sous-catégorique propre doit vérifier en outre
\[(Id_A)^\alpha = Id_{A^\alpha},\]
et pour une dynamique catégorique on doit avoir de plus
\[(f\circ g)^\alpha = f^\alpha \odot g^\alpha.\]  Pour les dynamorphismes entre dynamiques de même moteur, il n'y a pas de contrainte supplémentaire, et pour les dynamorphismes plus généraux de la forme
${(\Delta,\delta):\alpha\looparrowright\beta}$ on demande que $\Delta$ soit un foncteur.

Enfin, les dynamiques sous-catégoriques incluses dans une dynamique sous-catégoriques constituent un treillis complet avec l'union pour $\sup$.

\subsection{Dynamiques sous-catégoriques scandées}

\subsubsection{Définitions}

\begin{df}  Une \emph{dynamique sous-catégorique scandée} $\tau: \alpha\looparrowright h$ est un dynamorphisme sous-catégorique $\tau$ tel que $Gr(\tau)$ soit une dynamique graphique scandée. 
\end{df}

Autrement dit,
une {dynamique sous-catégorique scandée sur une petite caté\-gorie $\mathbf{C}$} est un dynamorphisme déterministe $\tau: \alpha\looparrowright h$ avec $\alpha$ une dynamique sous-catégorique sur $\mathbf{C}$ et $h$ une horloge sur $\mathbf{C}$. 

Comme dans le cas des dynamiques graphiques, la dynamique sous-caté\-gorique scandée $\tau$ est également appelée une \emph{datation} ou une \emph{scansion}, $\alpha$ est appelée la dynamique sous-catégorique de $\tau$ et $h$ est son horloge. S'il n'y a pas d'ambiguïté sur la datation, la dynamique sous-catégorique $\tau$ sera parfois désignée par sa dynamique sous-catégorique $\alpha$.

\subsubsection{Dynamorphismes scandés}

\begin{df}\label{df dynamo sc scandes}[Dynamorphismes sous-catégoriques scandés]
 Étant données deux dynamiques scandées $\rho$ et $\tau$, on appelle \emph{dynamorphisme sous-catégorique scandé}, ou simplement dynamorphisme, de $\rho$ vers $\tau$, tout dynamorphisme graphique scandé 
\[
(\Delta,\delta,d):Gr(\rho)\looparrowright Gr(\tau)
\]
tel que $\Delta$ soit un foncteur du moteur de $\rho$ vers le moteur de $\tau$. $\Delta$ est appelé la partie fonctorielle du dynamorphisme scandé $(\Delta,\delta,d)$, $\delta$ sa partie transitionnelle, et $d$ sa partie horloge.
\end{df}

Autrement dit, un dynamorphisme de 
\[
\rho:
(\alpha:\mathbf{C}\rightharpoondown \mathbf{P})
\looparrowright
(h:\mathbf{C}\rightarrow \mathbf{P})
\]
vers 
\[  
\tau:
(\beta:\mathbf{D}\rightharpoondown \mathbf{P})
\looparrowright
(k:\mathbf{D}\rightarrow \mathbf{P})
\] 
est un triplet $(\Delta,\delta,d)$ tel que 
\begin{enumerate}
\item $(\Delta,\delta)$ est un dynamorphisme sous-catégorique de $\alpha$ vers $\beta$,
\item $(\Delta,d)$ est un dynamorphisme de l'horloge $h$ vers $k$,
\item pour tout $S\in\dot{\mathbf{C}}$, la condition de synchronisation entre $\rho$ et $\tau$ est satisfaite :
\[\tau_{\Delta_S}\odot \delta_S\subset d_S\odot \rho_S. \]  
\end{enumerate}

\begin{rmq}\label{rmq conditions synchro et real dyna scand} Rappelons ce que nous avons indiqué dans les sections §2.2.2 et §2.5.2 de \cite{Dugowson:20150807}, et en particulier dans la remarque 7 du texte en question, à savoir  que la condition de synchronisation écrite ci-dessus permet d'interpréter les réalisations d'une dynamique scandée $\tau:\alpha\looparrowright h$ comme les dynamorphismes vers $\tau$ d'une certaine dynamique scandée $[h]$ canoniquement associée à $h$. 
\end{rmq}

\subsection{Multidynamiques sous-catégoriques}

\subsubsection{Définitions}

\begin{df}
Une \emph{$\mathbf{C}$-multi-dynamique sous-catégorique} est une famille $\alpha=(\alpha_\mu:\mathbf{C}\rightharpoondown \mathbf{P})_{\mu\in M}$ de dynamiques sous-catégoriques de même moteur $\mathbf{C}$, indexée par un ensemble non vide $M$ appelé ensemble des \emph{paramètres} de la multi-dynamique, telle que  $Gr(\alpha)$ soit une $Gr(\mathbf{C})$-multi-dynamique graphique, où l'on désigne par
$Gr(\alpha)$ la famille de dynamiques graphiques $Gr(\alpha_\mu)_{\mu\in M}$.
Nous dirons que la multi-dynamique sous-catégorique $\alpha$ est \emph{propre} (respectivement \emph{catégorique}) si \emph{pour tout} $\mu\in M$, la dynamique sous-catégorique $\alpha_\mu$ est propre (respectivement catégorique).
\end{df}

D'après la définition des multi-dynamiques graphiques\footnote{Voir \cite{Dugowson:20150807}, définition 12.}, $M$ étant un ensemble non vide, une famille $\alpha=(\alpha_\mu:\mathbf{C}\rightharpoondown \mathbf{P})_{\mu\in M}$ de dynamiques sous-catégoriques de même moteur $\mathbf{C}$ est une muti-dynamique sous-catégorique si et seulement si pour tout $S\dot{\mathbf{C}}$ il existe un ensemble $S^\alpha$ tel que
\[\forall \mu\in M, S^{\alpha_\mu}=S^\alpha.\] Bien entendu, l'ensemble
\[st(\alpha)=  \bigcup_{S\in\dot{\mathbf{C}}}{S^\alpha}\]
est appelé l'ensemble des états de la multi-dynamique $\alpha$.

 Conformément aux notations introduites dans \cite{Dugowson:20150807}, 
nous écrirons
$\alpha:\mathbf{C}\rightharpoondown \mathbf{P}^{\underrightarrow{{M}}}$ pour exprimer que $\alpha$ est une  multi-dynamique sous-catégorique de moteur $\mathbf{C}$ et d'ensemble de paramètre $M$. 
Pour toute flèche  $(e:S\rightarrow T)\in\overrightarrow{\mathbf{C}}$ et tout  paramètre $\mu\in M$, on notera indifféremment  $e^{\alpha_\mu}$ ou $e^\alpha_\mu$ la transition de paramètre $\mu$ associée par $\alpha$ à $e$, et  la famille de transitions $e^\alpha$ 
sera souvent désignée  par la notation
\[e^\alpha: S^\alpha \twosquig_M  T^\alpha,\] ou plus simplement, s'il n'y pas d'ambiguïté sur les paramètres, par
\[e^\alpha: S^\alpha \twosquig  T^\alpha.\]
Les deux expressions  
$(e^\alpha_\mu:S^\alpha\rightsquigarrow T^\alpha)_{\mu\in M}$ et $e^\alpha:S^\alpha \twosquig_M  T^\alpha$ signifierons donc toutes deux la même chose, à savoir que $e^\alpha$ est une famille indexée par $M$ de transitions de $S^\alpha$ vers $T^\alpha$.
 
Par opposition aux multi-dynamiques, les dynamiques sous-catégoriques don\-nées par la définition \ref{df dyscat} et qui s'identifient aux multi-dynamiques ayant un singleton pour ensemble de paramètres, seront parfois appelées des \textit{mono-dynamiques}.

Récapitulons. Une multi-dynamique sous-catégorique $\alpha:\mathbf{C}\rightharpoondown \mathbf{P}^{\underrightarrow{{M}}}$ consiste en la donnée 
\begin{itemize}
\item d'une application qui à tout objet $S\in\dot{\mathbf{C}}$ associe un ensemble $S^\alpha$, de telle sorte que $S\neq T \Rightarrow S^\alpha\cap T^\alpha=\emptyset$,
\item d'une application qui à toute flèche $(S\stackrel{f}{\rightarrow}T)\in\overrightarrow{\mathbf{C}}$ associe une famille de transitions $f^\alpha: S^\alpha \twosquig_M  T^\alpha$ indexée par $M$ de telle sorte que pour tout $S\in\dot{\mathbf{C}}$ on ait
\[
\forall \mu\in M, (Id_S)^\alpha_\mu\subset Id_{S^\alpha}
\]
et pour tout couple $(g,f)$ de flèches composables de $\mathbf{C}$, on ait
\[
\forall \mu\in M, (g\circ f)^\alpha_\mu\subset g^\alpha_\mu \odot f^\alpha_\mu.
\]
\end{itemize}

\begin{rmq}\label{rmq nettoyage multi-dynamiques sc} L'opération de nettoyage décrite en section §\ref{subsubsec Dyna SCat propres} qui, à partir d'une dynamique sous-catégorique produit une dynamique sous-catégorique propre par le retrait de tous les états \og hors du coup\fg\, ne fonctionne plus aussi bien pour les multi-dynamiques sous-catégoriques. En effet, un état peut fort bien être dans le coup pour certaines valeurs du paramètre et hors du coup pour d'autres valeurs, de sorte qu'en appelant à présent \emph{nettoyage} l'opération consistant à retirer tous les états qui sont hors du coup pour \emph{toutes} les valeurs du paramètre,  on obtient une multi-dynamique que l'on pourrait dire \emph{semi-propre}, à savoir telle que pour tout état il existe une valeur du paramètre pour laquelle cet état soit dans le coup, mais qui en général n'est pas propre au sens propre.  
\end{rmq}

\begin{df}\label{df determinisme multi-dynamique sc} Une multi-dynamique sous-catégorique 
$\alpha:\mathbf{C}\rightharpoondown \mathbf{P}^{\underrightarrow{{M}}}$ est dite déterministe (respectivement quasi-déterministe) si pour tout $\mu\in M$, la mono-dynamique $\alpha_\mu$ est déterministe (repsectivement quasi-déterministe).
\end{df}

D'après la proposition \ref{prop dysc deterministe implique categorique}, une multi-dynamique sous-catégorique déterministe est nécessairement catégorique.

\subsubsection{Multi-dynamorphismes sous-catégoriques}

La catégorie des multi-dynamiques sous-catégoriques a pour objets toutes les multi-dynamiques sous-catégoriques $\alpha:\mathbf{C}\rightharpoondown \mathbf{P}^{\underrightarrow{{M}}}$,
 où $\mathbf{C}$ décrit la classe des petites catégories et $M$ celle des ensembles, et pour flèches les multi-dyna\-morphismes sous-catégoriques définis de la façon suivante.

\begin{df}\label{df multi-dynamorphismes sc}
Étant données $\alpha:\mathbf{C}\rightharpoondown\mathbf{P}^{\underrightarrow{{L}}}$ et
$\beta:\mathbf{D}\rightharpoondown \mathbf{P}^{\underrightarrow{M}}$ deux multi-dynamiques sous-catégoriques, un \emph{multi-dynamorphisme sous-catégorique} --- ou plus simplement un dynamorphisme ---  de $\alpha$ vers $\beta$ est un triplet $(\theta,\Delta,\delta)$  qui soit un dynamorphisme (au sens des multi-dynamiques graphiques) de $Gr(\alpha)$ vers $Gr(\beta)$ et tel que $\Delta$ soit un foncteur $\mathbf{C}\rightarrow \mathbf{D}$.
\end{df}

Autrement dit, un dynamorphisme 
\[
(\theta,\Delta,\delta):
(\alpha:\mathbf{C}\rightharpoondown\mathbf{P}^{\underrightarrow{{L}}})
\looparrowright
(\beta:\mathbf{D}\rightharpoondown \mathbf{P}^{\underrightarrow{M}})
\]
est constitué
\begin{itemize}
\item d'une application $\theta:L\rightarrow M$,
\item d'un foncteur $\Delta:\mathbf{C}\rightarrow\mathbf{D}$,
\item d'une transition 
$\delta:st(\alpha)\rightsquigarrow st(\beta)$,
\end{itemize} 
\noindent telles que, pour tout $\lambda\in L$, $(\Delta,\delta)$ définit un dynamorphisme de la dynamique sous-catégorique $\alpha_\lambda$ vers la dynamique sous-catégorique $\beta_{\theta(\lambda)}$.

Ainsi, pour tout $\lambda\in L$, tous $S$ et $T$ dans $\dot{\mathbf{C}}$ et tout $(f:S\rightarrow T)\in \overrightarrow{\mathbf{C}}$, on doit avoir
\[\delta_T\odot f^\alpha_\lambda\subset (\Delta f)^\beta_{\theta(\lambda)}\odot\delta_S.\] 

\begin{rmq}\label{rmq C-multi-dyna sc vers mono} [$\mathbf{C}$-multi-dynamorphismes]
Dans le cas particulier où $\mathbf{C}=\mathbf{D}$, on entendra implicitement et sauf mention contraire par \emph{dynamorphisme} 
\[
(\alpha:\mathbf{C}\rightharpoondown\mathbf{P}^{\underrightarrow{{L}}})
\looparrowright
(\beta:\mathbf{C}\rightharpoondown \mathbf{P}^{\underrightarrow{M}})
\]
un $\mathbf{C}$-dynamorphisme, autrement dit un multi-dynamorphisme sous-catégorique $(\theta,\Delta,\delta)$ avec pour $\Delta$ le foncteur identité
\[\Delta=Id_\mathbf{C}:\mathbf{C}\rightarrow\mathbf{C}.\] 

Dans ces conditions, le lemme suivante découle immédiatement de la définition \ref{df multi-dynamorphismes sc}.

\begin{lm}\label{lm multi-dynamorphisme} Étant données $\mathbf{C}$  une petite catégorie, et $a:Gr(\mathbf{C})\rightarrow Gr(\mathbf{P}^{\underrightarrow{{L}}})$ et $b:Gr(\mathbf{C})\rightarrow Gr(\mathbf{P}^{\underrightarrow{{M}}})$ deux multi-dynamiques graphiques sur  $Gr(\mathbf{C})$,  un multi-dynamorphisme graphique $(\theta,Id_{Gr(\mathbf{C})},\delta):a\looparrowright b$  est sous-catégorique si et seulement si $a$ et $b$ sont sous-catégoriques.
\end{lm}

\paragraph{$\mathbf{C}$-multi-dynamorphismes vers une mono-dynamique.} En particulier, comme ce sera le cas ci-après dans la section §\ref{subsubsec df DySCO}, si la dynamique d'arrivée $\beta$ est une monodynamique de même moteur que $\alpha$, l'application $\theta$ est nécessairement l'unique application $L\rightarrow \{*\}$ de sorte qu'un dynamorphisme sous-catégorique $\alpha\looparrowright\beta$ se réduit 
à la donnée de la partie transitionnelle d'un dynamorphisme graphique
\[\delta:
(\alpha:\mathbf{C}\rightharpoondown\mathbf{P}^{\underrightarrow{{L}}})
\looparrowright
(\beta:\mathbf{C}\rightharpoondown \mathbf{P})
\] vérifiant à ce titre 
\[\forall \lambda\in L, \forall (f:S\rightarrow T)\in \overrightarrow{\mathbf{C}}, \delta_T\odot f^\alpha_\lambda\subset f^\beta\odot\delta_S,\] où comme d'habitude $\delta_U$ désigne la restriction selon $U^\alpha\rightarrow U^\beta$ de l'application $\delta:st(\alpha)\rightarrow st(\beta)$.
\end{rmq}

\subsubsection{Quotients paramétriques}\label{subsubsec reduc param mdsc}

Le quotient paramétrique $\alpha/{\sim}$ d'une multi-dynamique sous-catégorique  $\alpha:\mathbf{C} \rightharpoondown \mathbf{P}^{\underrightarrow{{M}}}$ par une relation d'équivalence $\sim$  sur $M$ est définie comme étant 
le quotient paramétrique de la multi-dynamique graphique\footnote{Voir la définition 14 dans \cite{Dugowson:20150807}.} $Gr(\alpha)$ par $\sim$. Ce quotient constitue en effet, en conséquence de la proposition \ref{prop union dyna sous-cat} portant sur l'union d'une famille de dynamiques sous-catégoriques, une dynamique sous-catégorique. Posant  $\beta=\alpha/{\sim}$ et $\widetilde{M}=M/{\sim}$, on a :
\begin{itemize}
\item pour tout sommet $S$ de $\mathbf{C}$, $S^\beta=S^\alpha$,
\item pour toute flèche $(f:S\rightarrow T)\in\overrightarrow{\mathbf{C}}$, pour toute classe $\lambda\in \widetilde{M}$ et tout état $a\in S^\beta$, 
\[f^\beta_\lambda(a)=\bigcup_{\mu\in\lambda}f^\alpha_\mu(a).\]
\end{itemize}

\subsection{Dynamiques sous-catégoriques ouvertes}

\subsubsection{Définitions}\label{subsubsec df DySCO}

On appelle \emph{dynamique sous-catégorique ouverte} (ou \emph{dynamique ouverte sous-catégorique}, ou simplement \emph{dynamique ouverte})  toute  multi-dynamique sous-catégorique scandée, autrement dit toute multi-dynamique sous-catégorique munie d'une horloge et d'une datation. Plus précisément, posons la définition suivante.

\begin{df}\label{df DySCO}  Une \emph{dynamique (sous-catégorique) ouverte} de moteur $\mathbf{C}$ est un $\mathbf{C}$-dynamorphisme 
\[
\tau:
(\alpha=(\alpha_\mu)_{\mu \in M}:\mathbf{C}\rightharpoondown \mathbf{P}^{\underrightarrow{{M}}} )
\looparrowright 
(h:\mathbf{C}\rightarrow \mathbf{P})\]
d'une $\mathbf{C}$-multi-dynamique $\alpha$ vers une $\mathbf{C}$-horloge $h$. Une telle dynamique ouverte est dite \emph{propre} (respectivement \emph{catégorique}) si $\alpha$ est une multi-dynamique sous-catégorique propre (respectivement catégorique).
\end{df}

Autrement dit, une dynamique sous-catégorique ouverte est un multi-dyna\-morphisme sous-catégorique $\tau$ tel que $Gr(\tau)$ soit une dynamique graphique ouverte.

De la remarque \ref{rmq C-multi-dyna sc vers mono}, et en particulier du lemme \ref{lm multi-dynamorphisme}, on déduit immédiatement le lemme suivant :
\begin{lm}\label{lm dygrouv vers dysco} Étant donnés $\mathbf{C}$ une petite catégorie et 
\[
\tau:
(a:Gr(\mathbf{C})\rightarrow Gr(\mathbf{P}^{\underrightarrow{{L}}}))
\looparrowright
(h:Gr(\mathbf{C})\rightarrow Gr(\mathbf{P}))
\] 
une dynamique graphique ouverte sur le graphe $Gr(\mathbf{C})$, une condition nécessaire et suffisante pour que $\tau:a\looparrowright h$ constitue une dynamique sous-catégorique ouverte sur $\mathbf{C}$ est que $a$ et $h$ soient sous-catégoriques sur $\mathbf{C}$.
\end{lm}

Par rapport à la notion de dynamique graphique ouverte, il n'y a pas donc pas de condition supplémentaire portant sur $\tau$. Conformément à la section §2.4.1 de \cite{Dugowson:20150807}, $\tau$ peut donc être vue comme une application $st(\alpha)\rightarrow st(h)$ astreinte à vérifier la condition suivante
\[
\forall \mu\in M,
\forall (f:S\rightarrow T)\in\overrightarrow{\mathbf{C}},
\forall a\in S^\alpha,
\forall b\in f^\alpha_\mu(a),
\tau(b)=f^h (\tau (a)).
\]

On désigne souvent la dynamique ouverte $A=(\tau:\alpha \looparrowright h)$ par sa multi-dynamique $\alpha$. Ainsi, l'ensemble des états de $A$, noté $st(A)$, n'est rien d'autre que $st(\alpha)$. De plus, $\tau$ est aussi appelé la scansion ou la datation de la dynamique ouverte $A$, $h$ l'horloge de $A$ et $M$ l'ensemble des paramètres de $A$. 

\begin{df}\label{df determinisme dyna sc ouverte} Une dynamique sous-catégorique ouverte
\[A=(
\tau:
(\alpha=(\alpha_\mu)_{\mu \in M}:\mathbf{C}\rightharpoondown \mathbf{P}^{\underrightarrow{{M}}} )
\looparrowright 
(h:\mathbf{C}\rightarrow \mathbf{P}))\] est dite \emph{déterministe} (respectivement \emph{quasi-déterministe}) si sa multi-dynamique $\alpha$ est déterministe (respectivement quasi-déterministe).
\end{df}

D'après la définition \ref{df determinisme multi-dynamique sc}, une dynamique sous-catégorique ouverte est donc déterministe lorsque chacune des mono-dynamiques $\alpha_\mu$ est déterministe. Et d'après la proposition \ref{prop dysc deterministe implique categorique}, une dynamique sous-catégorique ouverte déterministe est donc nécessairement catégorique.

\subsubsection{Dynamorphismes de dynamiques ouvertes sous-catégoriques}
\label{subsubsec DySCO}

On constitue la catégorie $\mathbf{DySCO}$ des dynamiques sous-catégoriques ouvertes en prenant pour flèches les multi-dynamorphismes scandés ainsi définis :

\begin{df}[Dynamorphismes ouverts]\label{df dyna DySCO} On appelle \emph{dynamorphisme ouvert} ou \emph{multi-dynamorphisme scandé}, ou plus simplement  \emph{dynamorphisme}, d'une dynamique sous-catégorique ouverte
\[
A=
(
\rho:
(\alpha:\mathbf{C}\rightharpoondown \mathbf{P}^{\underrightarrow{L}})
\looparrowright 
(h:\mathbf{C}\rightarrow \mathbf{P})
)
\]
vers une dynamique sous-catégorique ouverte
\[
B=
(
\tau:
(\beta:\mathbf{D}\rightharpoondown \mathbf{P}^{\underrightarrow{M}})
\looparrowright 
(k:\mathbf{D}\rightarrow \mathbf{P})
)
\]
la donnée d'un quadruplet $(\theta,\Delta,\delta,d)$ tel que
\begin{enumerate}
\item $(\theta,\Delta,\delta)$ est un multi-dynamorphisme sous-catégorique de $\alpha$ vers $\beta$,
\item $(\Delta,d)$ est un dynamorphisme de $h$ vers $k$,
\item pour tout $S\in\dot{\mathbf{C}}$, la condition suivante de synchronisation entre $\rho$ et $\tau$ est satisfaite :
\[\tau_{\Delta_S}\odot \delta_S\subset d_S\odot \rho_S. \]  
\end{enumerate}
\end{df}

Autrement dit, un dynamorphisme $(\theta,\Delta,\delta,d):A\looparrowright B$
est un quadruplet $(\theta,\Delta,\delta,d)$ tel que $Gr(\theta,\Delta,\delta,d)= (\theta, Gr(\Delta), \delta,d)$ soit un dynamorphisme de la dynamique graphique ouverte $Gr(A)$ vers la dynamique graphique ouverte $Gr(B)$, celles-ci étant simplement obtenues en oubliant la possibilité de composer les écoulements de $\mathbf{C}$ et $\mathbf{D}$. En particulier, conformément à la définition 16 de \cite{Dugowson:20150807}, étant donné un dynamorphisme sous-catégorique ouvert $(\theta,\Delta,\delta,d)$, on appellera
\begin{itemize}
\item $\theta$ sa  \emph{partie paramétrique},
\item $\Delta$ sa  \emph{partie fonctorielle},
\item $\delta$ sa  \emph{partie transitionnelle},
\item et $d$ sa \emph{partie horloge}.
\end{itemize}

Les mono-dynamiques pouvant être considérées comme des multi-dynamiques particulières et  toute dynamique pouvant être canoniquement scandée par l'horloge essentielle\footnote{Voir l'exemple 4 de \cite{Dugowson:20150807}.}, on vérifie immédiatement que la définition \ref{df dyna DySCO} ci-dessus généra\-lise toutes les définitions de dynamorphismes sous-catégoriques données précé\-demment.

\subsubsection{Quotient paramétrique d'une dynamique ouverte}
\label{subsubsec reduc param mdsco}

Grâce à la notion de quotient paramétrique d'une multi-dynamique sous-catégorique présentée en section §\ref{subsubsec reduc param mdsc}, elle-même rendue possible par la proposition \ref{prop union dyna sous-cat}, on  définit ainsi le quotient paramétrique d'une dynamique sous-catégorique ouverte par une relation d'équivalence sur l'ensemble de ses paramètres :

\begin{df}\label{df reduc param ouverte sc}
Étant donnée $\tau:(\alpha:\mathbf{C}\rightharpoondown \mathbf{P}^{\underrightarrow{{M}}})\looparrowright h)$ une dynamique ouverte sous-catégorique de moteur $\mathbf{C}$, d'horloge $h$, de datation $\tau$, d'ensemble de paramètres $M$ et de multi-dynamique $\alpha$, et $\sim$ une relation d'équivalence sur $M$, on appelle quotient de $\alpha$ par $\sim$ et l'on note  $\alpha/{\sim}$ la dynamique sous-catégorique ouverte ayant le même moteur, la même horloge, les mêmes états et la  même datation que $\alpha$, dont l'ensemble des paramètres est $\widetilde{M}=M/{\sim}$  et dont la multi-dynamique est $\alpha/{\sim}$.
\end{df}

\subsection{Réalisations}\label{subsec realisations sc}

Dans la présente section §\ref{subsec realisations sc}, la notion de \emph{réalisation} vue dans le cas des dynamiques graphiques (section §2.5 de \cite{Dugowson:20150807}) et dans celui des  dynamiques sous-catégoriques (section §\ref{subsubsec horloge cat success real} du présent document) est élargie  à la version sous-catégorique des multi-dynamiques, des dynamiques scandées et des dynamiques ouvertes.


\subsubsection{Cas des dynamiques sous-catégoriques}

Conformément aux notations de la section §2.5.1. de \cite{Dugowson:20150807}, nous noterons $\mathcal{S}_{(h,\alpha)}$ l'ensemble des $h$-réalisations d'une $\mathbf{C}$-dynamique sous-catégorique $\alpha$, où $h$ est une $\mathbf{C}$-horloge. Rappelons\footnote{Voir la  définition \ref{df realisation de DySC pour horloge}.} qu'une telle réalisation est un  dynamorphisme quasi-déterministe de $h$ dans $\alpha$. Concernant les conventions d'écriture relatives aux transitions quasi-déterministes, voir plus haut la note de bas de page \ref{foot ecriture transition quasidet}, ainsi que, dans \cite{Dugowson:20150807}, la définition 5 et la section §2.5.1.

\subsubsection{Cas des dynamiques scandées}

\begin{df} On appelle \emph{réalisation d'une dynamique sous-catégorique scandée} $A$ toute réalisation de la dynamique graphique scandée $Gr(A)$.
\end{df}

Autrement dit, si $A=(\tau:\alpha \looparrowright h)$ est une dynamique sous-catégorique scandée, une réalisation de $A$ est une réalisation de la dynamique graphique scandée $Gr(\tau):Gr(\alpha) \looparrowright Gr(h)$, à savoir 
 une fonction $\mathfrak{a}:st(h)\dashrightarrow st(\alpha)$ définie sur une partie $df(\mathfrak{a})\subset st(h)$ et qui vérifie les deux propriétés suivantes :
\[\forall t\in df(\mathfrak{a}), \tau(\mathfrak{a}(t))=t,\]
et
\[\forall (S\stackrel{f}{\rightarrow}T)\in\overrightarrow{\mathbf{C}}, 
 \forall t\in S^h, 
 f^h(t)\in df(\mathfrak{a})\Rightarrow 
 t\in df(\mathfrak{a}) 
\,\mathrm{et}\,
\mathfrak{a}(f^h(t))\in f^\alpha (\mathfrak{a}(t)),\]
cette dernière propriété impliquant notamment celle-ci:
\[\forall (S\stackrel{f}{\rightarrow}T)\in\overrightarrow{\mathbf{C}}, \forall t\in S^h, \, t\notin df(\mathfrak{a})\Rightarrow f^h(t)\notin df(\mathfrak{a}).
\]

On note $\mathcal{S}_A$ l'ensemble des réalisations de $A$.

\begin{rmq} Comme rappelé plus haut dans la remarque \ref{rmq conditions synchro et real dyna scand} en référence aux sections §2.2.2 et §2.5.2 de \cite{Dugowson:20150807} et en particulier à la remarque 7 du texte en question, une telle réalisation d'une dynamique scandée peut être vue comme un dynamorphisme scandé de $[h]$ vers $A$, où $[h]$ est une dynamique scandée canoniquement associée à $h$, et cela grâce à la condition de synchronisation figurant dans la définition \ref{df dynamo sc scandes}, condition qui se trouve ainsi justifiée.
\end{rmq}

\subsubsection{Cas des multi-dynamiques}
\label{subsubsec h realisation de multi-dynamique sc}

\begin{df} Étant donnée une $\mathbf{C}$-horloge $h$, on appelle \emph{$h$-réalisation}  d'une $\mathbf{C}$-multi-dynamique sous-catégorique $\alpha:\mathbf{C}\rightharpoondown \mathbf{P}^{\underrightarrow{{L}}}$ toute $Gr(h)$ réalisation de la $Gr(\mathbf{C})$-multi-dynamique graphique $Gr(\alpha)$.
\end{df}

Autrement dit, une $h$-réalisation de $\alpha$ est un  $\mathbf{C}$-multi-dynamorphisme quasi-déterministe de la mono-dynamique $h$ dans $\alpha$, ce qui revient à dire qu'elle consiste la donnée d'une valeur du paramètre $\lambda\in L$ et d'une $h$-réalisation de la mono-dynamique sous-catégorique $\alpha_\lambda$. 
On peut donc  aussi plus simplement voir  une telle réalisation comme un couple $(\lambda,\mathfrak{a})$ constitué d'une valeur  $\lambda\in L$ et 
d'une fonction $\mathfrak{a}:st(h)\supset df(\mathfrak{a})\dashrightarrow st(\alpha)$ vérifiant
\[\forall (S\stackrel{f}{\rightarrow}T)\in\overrightarrow{\mathbf{C}}, 
\forall t\in S^h, 
f^h(t)\in df(\mathfrak{a})\Rightarrow 
t\in df(\mathfrak{a}) 
\,\mathrm{et}\,
\mathfrak{a}(f^h(t))\in f^\alpha_\lambda (\mathfrak{a}(t)).\]

Étant donnée $(\lambda,\mathfrak{a})$ une $h$-réalisation de la multi-dynamique sous-catégorique $\alpha$, nous appellerons $\lambda$ la \emph{partie interne} ou \emph{paramétrique} de cette réalisation, tandis que $\mathfrak{a}$ sera appelée sa partie externe. 

Conservant les notations introduites dans la section §2.5.3 de \cite{Dugowson:20150807}, nous noterons ${\mathcal{S}}_{(h,\alpha_\lambda)}$ l'ensemble des $h$-réalisations de la mono-dynamique sous-catégorique $\alpha_\lambda$, où $\lambda$ est une valeur donnée du paramètre dont dépend la multi-dynamique $\alpha$, et
nous noterons
${\mathcal{S}}_{(h,\alpha)}$
l'ensemble des \emph{parties externes} des $h$-réalisations de la $\mathbf{C}$-multi-dynamique sous-catégorique $\alpha=(\alpha_\lambda)_{\lambda\in L}$, de sorte que l'on a
\[
{\mathcal{S}}_{(h,\alpha)}
=\bigcup_{\lambda\in L} {\mathcal{S}}_{(h,\alpha_\lambda)}.\]
Dans la formule ci-dessus, l'union \emph{n'est pas disjointe}. Par exemple, la réalisation vide est commune à tous les ${\mathcal{S}}_{(h,\alpha_\lambda)}$.

\subsubsection{Cas des dynamiques ouvertes}
\label{subsubsec real dyna ouv sc}

\begin{df}\label{df realisation dysc ouverte} Si $A$ est une dynamique sous-catégorique ouverte, on appelle \emph{réalisation de $A$} toute réalisation de la dynamique graphique ouverte $Gr(A)$.
\end{df}

Autrement dit, étant donnée
\[
A=(
\tau:
(\alpha=(\alpha_\lambda)_{\lambda \in L}:\mathbf{C}\rightharpoondown \mathbf{P}^{\underrightarrow{{L}}} )
\looparrowright 
(h:\mathbf{C}\rightarrow \mathbf{P})
)
\]
une $\mathbf{C}$-dynamique ouverte, une
réalisation de $A$ consiste en un couple $(\lambda,\mathfrak{a})$ constitué d'une valeur $\lambda\in L$ et d'une fonction $\mathfrak{a}:st(h)\dashrightarrow st(\alpha)$ définie sur une partie $df(\mathfrak{a})\subset st(h)$ et qui vérifie les deux propriétés suivantes :
\[\forall t\in df(\mathfrak{a}), \tau(\mathfrak{a}(t))=t,\]
et
\[\forall (S\stackrel{f}{\rightarrow}T)\in\overrightarrow{\mathbf{C}}, 
 \forall t\in S^h, 
 f^h(t)\in df(\mathfrak{a})\Rightarrow 
 t\in df(\mathfrak{a}) 
\,\mathrm{et}\,
\mathfrak{a}(f^h(t))\in f^\alpha_\lambda (\mathfrak{a}(t)).\]

Comme indiqué ci-dessus en section §\ref{subsubsec h realisation de multi-dynamique sc} dans le cas des multi-dynamiques sous-catégoriques, $\lambda$ sera appelé la \emph{partie interne} ou \emph{paramétrique} de la réalisation $(\lambda,\mathfrak{a})$ de $A$, tandis que $\mathfrak{a}$ est la \emph{partie externe} de cette réalisation.

Nous noterons $\mathcal{S}_A$ l'ensemble des parties externes des réalisations de la dynamique sous-catégorique ouverte $A$, et $\mathcal{S}_{(A,\lambda)}$ ou $\mathcal{S}_{A_\lambda}$ l'ensemble des réalisations de la mono-dynamique scandée $A_\lambda=(
\tau:
(\alpha_\lambda:\mathbf{C}\rightharpoondown \mathbf{P})
\looparrowright 
(h:\mathbf{C}\rightarrow \mathbf{P})
)$, de sorte que 
\[\mathcal{S}_A=\bigcup_{\lambda\in L}\mathcal{S}_{A_\lambda},\]
cette union n'étant pas en général disjointe.

\begin{rmq}\label{rmq usage realisation au lieu de partie externe} Souvent, et sans que  cela ne porte  à conséquence, nous parlerons de \emph{la réalisation $\mathfrak{a}$ de $A$} au lieu de 
\emph{la partie externe $\mathfrak{a}$ d'une réalisation de $A$}. Cette façon de parler conduira par exemple à désigner l'ensemble $\mathcal{S}_A$ comme \og ensemble des réalisations de $A$\fg\, bien que l'expression soit impropre et qu'en l'occurrence il vaille mieux l'éviter.
\end{rmq}

\subsubsection{Réalisations passant par un état}
\label{subsubsec realisations sc passant par etat}

\begin{df}
Étant donnée une dynamique sous-catégorique ouverte $A$,  nous dirons qu'une réalisation\footnote{Voir la remarque \ref{rmq usage realisation au lieu de partie externe} ci-dessus.} $\mathfrak{a}$ de $A$ passe par un état $a\in st(A)$, si elle passe par $a$ en tant que réalisation de la dynamique graphique $Gr(A)$. 
\end{df}

Comme pour les dynamiques graphiques, nous écrirons
\[\mathfrak{a}\rhd a,\] 
pour exprimer que $\mathfrak{a}$ passe par $a$.
 
 Ainsi, pour $A=(\tau:(\alpha_\lambda)_{\lambda\in L}\looparrowright h)$, on a 
\[\mathfrak{a}\rhd a \Leftrightarrow \mathfrak{a}(\tau(a))=a.\]

Dans le cas où $a$ et $b$ sont deux états de la dynamique ouverte $A$ tels que $\tau(b)$ succède\footnote{Voir la section §\ref{subsubsec horloge cat success real}.} à $\tau(a)$, nous écrirons
\[\mathfrak{a}\rhd a, b\] pour exprimer que $\mathfrak{a}$ passe par $a$ \emph{puis} qu'elle passe par $b$.

\section{Interactions et familles dynamiques sous-caté\-goriques}
\label{section familles dynamiques sc}

La définition d'une interaction sur une famille de dynamiques ouvertes et la définition d'une famille dynamique ont été données dans \cite{Dugowson:20150807}, section §3, dans le cadre des dynamiques graphiques. Nous les reprenons ici dans le cadre des dynamiques sous-catégoriques, quasiment à l'identique. Ces définitions s'appuient notamment sur la notion de \emph{relation binaire multiple}, introduite dans la section §1 de \cite{Dugowson:20150807} et qui s'appuie elle-même sur les \emph{relations multiples} considérées dans notre texte \cite{Dugowson:201505}, où se trouve en particulier définie la structure connective d'une telle relation multiple.

\begin{rmq} Soulignons que, par rapport à notre conférence \cite{Dugowson:20130521} du 21 mai 2013, les interactions considérées ici sont considérablement plus générales puisqu'elles concernent de façon globale toute famille, finie ou non, de dynamiques, tandis que dans la conférence en question nous ne considérions que l'influence d'une dynamique sur une autre.
\end{rmq}

Les lettres $A$, $B$, etc. désignerons donc des dynamiques sous-catégoriques ouvertes :
\[
A=
(
\tau:
(\alpha:\mathbf{C}\rightharpoondown \mathbf{P}^{\underrightarrow{L}})
\looparrowright 
(h:\mathbf{C}\rightarrow \mathbf{P})
),
\]
\[
B=
(
\rho:
(\beta:\mathbf{D}\rightharpoondown \mathbf{P}^{\underrightarrow{M}})
\looparrowright 
(k:\mathbf{D}\rightarrow \mathbf{P})
), \mathrm{etc...}
\]

Dans la définition \ref{df interaction sc} ci-dessous, $(A_i)_{i\in I}$ désigne une famille indexée par un ensemble $I$ de dynamiques sous-catégoriques ouvertes, avec
\[
A_i=
(
\tau_i:
(\alpha_i:\mathbf{C}_i\rightharpoondown \mathbf{P}^{\underrightarrow{L_i}})
\looparrowright 
(h_i:\mathbf{C}_i\rightarrow \mathbf{P})
),
\]
et, pour tout $i\in I$, $\mathcal{S}_{A_i}$ désigne, conformément à la section §\ref{subsubsec real dyna ouv sc}, l'ensemble des parties externes des réalisations de la dynamique sous-catégorique ouverte $A_i$.

\begin{df}\label{df interaction sc} $(A_i)_{i\in I}$ désignant comme ci-dessus une famille indexée par un ensemble $I$ de dynamiques sous-catégoriques ouvertes, on appelle \emph{interaction} pour cette famille la donnée d'une relation binaire multiple\footnote{Voir \cite{Dugowson:20150807}, section §1.2.} non vide $R\in \mathcal{BM}_{(\mathcal{S},\mathcal{L})}$, avec
\[\mathcal{S}=(\mathcal{S}_{A_i})_{i\in I} 
\quad\mathrm{et}\quad
\mathcal{L}=(L_i)_{i\in I},
\] vérifiant la propriété de cohérence suivante :
\[\forall (\mathfrak{a}_i,\lambda_i)_{i\in I}\in R,\,
 \forall i\in I,\,
 \mathfrak{a}_i\in \mathcal{S}_{(A_i,\lambda_i)}.\]
\end{df}

\pagebreak[3]

La définition suivante précise ce que nous entendrons par \emph{famille dynamique}, à savoir non seulement la donnée d'une famille de dynamiques ouvertes sous-catégoriques, mais aussi celle d'une interaction entre elles et d'une synchronisation globale assurée par l'une des dynamiques en jeu jouant en quelque sorte le rôle de chef d'orchestre.

\begin{df}\label{df famille dynamique sc} Une \emph{famille dynamique sous-catégorique}   $\mathcal{F}$ consiste en la donnée d'un quintuplet
\[\mathcal{F}=(I,i_0, (A_i)_{i\in I}, R, (\Delta_i,\delta_i)_{i\neq i_0})\] constitué
\begin{itemize}
\item d'un ensemble $I$ non vide, appelé \emph{index} de $\mathcal{F}$,
\item d'un élément $i_0\in I$, appelé \emph{indice synchronisateur} de $\mathcal{F}$,
\item d'une famille indexée par $I$ de dynamiques sous-catégoriques ouvertes $(A_i)_{i\in I}$  appelées \emph{composantes}  de $\mathcal{F}$, 
\item d'une interaction  $R\in \mathcal{BM}_{(\mathcal{S},\mathcal{L})}$ pour la famille $(A_i)_{i\in I}$, appelée \emph{interaction} de $\mathcal{F}$,
\item d'une famille 
\[
((\Delta_i,\delta_i) 
: h_{i_0}
\looparrowright 
h_i)_{i\in I\setminus\{i_0\}}
\] 
de dynamorphismes catégoriques\footnote{Il s'agit de dynamorphismes sous-catégoriques au sens de la définition \ref{df dynamorphismes dyscat}, mais les horloges étant des dynamiques catégoriques, nous pouvons  qualifier ces dynamorphismes eux-mêmes de catégoriques.} déterministes, appelés \emph{synchronisations} de  $\mathcal{F}$.
\end{itemize}
La composante $A_{i_0}$ de $\mathcal{F}$ sera appelée sa \emph{composante synchronisatrice}, et le moteur de cette composante sera appelé son \emph{moteur synchronisateur}\footnote{Dans les notations précédentes, le moteur synchronisateur de la famille $\mathcal{F}$ est donc la petite catégorie $\mathbf{C}_{i_0}$.}.
\end{df}

Dans la suite, lorsque nous parlerons de dynamiques ouvertes et de familles dynamiques, il s'agira toujours, sauf mention contraire, de dynamiques ouvertes sous-catégoriques et de familles dynamiques sous-catégoriques.

Grâce à la notion de structure connective d'une relation binaire multiple\footnote{Voir \cite{Dugowson:20150807}, section §1.2.3.}, nous pouvons finalement définir la structure connective d'une famille dynamique.

\begin{df} La \emph{structure connective d'une famille dynamique} est la structure connective de sa relation binaire multiple. L'\emph{ordre
connectif d'une famille dynamique} est l'ordre connectif\footnote{Sur la notion d'ordre connectif, voir par exemple, dans le cas fini, la définition 16 de \cite{Dugowson:201012}, et dans le cas général, la section §1.11 de \cite{Dugowson:201112} et \cite{Dugowson:201203}.} de sa structure connective.
\end{df}

\section{Dynamiques ouvertes engendrées par une famille dynamique $\mathcal{F}$}
\label{section DySCO engendrees}

\subsection{Rappels et notations} Le début de la section §4 de \cite{Dugowson:20150807} rappelle notamment qu'à toute relation binaire multiple  $R\in \mathcal{BM}_{(\mathcal{S}, \mathcal{L})}$ avec  $\mathcal{S}=(\mathcal{S}_{A_i})_{i\in I}$ et $\mathcal{L}= (L_i)_{i\in I}$ se trouve associée une relation binaire $rb(R):\Pi_I(\mathcal{S}_{A_i})\rightsquigarrow \Pi_I(L_i)$ dont l'image $\mathrm{Im}(rb(R))$ est définie par
\[\mathrm{Im}(rb(R))
=
\{
(\lambda_i)_{i\in I}\in \Pi_I(L_i),
\exists (\mathfrak{a}_i)_{i\in I}\in \Pi_I(\mathcal{S}_{A_i}),
(\mathfrak{a}_i,\lambda_i)_{i\in I}\in R
\}
,\] et que, pour tout $\mu\in\Pi_I(L_i)$, l'expression $rb(R)^{-1}(\mu)$ désigne l'ensemble
\[rb(R)^{-1}(\mu)=\{\mathfrak{a}\in\Pi_I(\mathcal{S}_{A_i}), (\mathfrak{a},\mu)\in rb(R)\},\] à la suite de quoi la définition 26 du texte en question associe à toute famille dynamique graphique, notons-la $\mathcal{G}$, une dynamique graphique ouverte  $[{\mathcal{G}}]_{\mathrm{p}}$ appelée la dynamique (graphique) ouverte  primo-engendrée par la famille $\mathcal{G}$.

\subsection{Dynamique graphique primo-engendrée par $\mathcal{F}$}\label{subsec dyna graph primo-engendree par famille dysc}

Soit maintenant $\mathcal{F}$ une famille dynamique sous-catégorique\footnote{Pour simplifier les notations, nous supposons que l'ensemble index $I$ contient un élément noté $0$ qui est choisi comme indice synchronisateur de la famille $\mathcal{F}$.}
\[\mathcal{F}=(I,0, (A_i)_{i\in I}, R, (\Delta_i,\delta_i)_{i\in I\setminus\{0\}}),\]  
avec, pour tout $i\in I$,
\[
A_i=
(
\tau_i:
(\alpha_i:\mathbf{C}_i\rightharpoondown \mathbf{P}^{\underrightarrow{L_i}})
\looparrowright 
(h_i:\mathbf{C}_i\rightarrow \mathbf{P})
).
\]

 Par les foncteurs d'oubli $Gr$ appliqués aux dynamiques sous-catégoriques ouvertes $A_i$ et aux foncteurs $\Delta_i$ (qui deviennent alors des \og foncteurs de graphes\fg), et tenant compte du fait --- sur lequel repose leur définition même\footnote{Voir plus haut la définition \ref{df realisation dysc ouverte}.}  ---  que les réalisations d'une dynamique sous-catégorique ouverte sont exactement les mêmes que celles de cette dynamique vue comme dynamique graphique ouverte, de sorte que l'interaction $R$ peut être conservée telle quelle dans le passage du point de vue \og sous-catégorique\fg\, au point de vue \og graphique\fg, on voit que $\mathcal{F}$ constitue une famille dynamique graphique, qu'en tant que telle nous noterons $Gr(\mathcal{F})$. D'où l'obtention d'une dynamique graphique ouverte $[Gr(\mathcal{F})]_{\mathrm{p}}$ : la dynamique graphique ouverte primo-engendrée par la famille dynamique sous-catégorique $\mathcal{F}$.
Selon la définition 26 de \cite{Dugowson:20150807}, on a 
\[[Gr(\mathcal{F})]_{\mathrm{p}}=({\rho:(\beta:\mathbf{G}\longrightarrow Gr(\mathbf{P}^{\underrightarrow{M}}))}\looparrowright k)\]
avec :
\begin{itemize}
\item $\mathbf{G}=Gr(\mathbf{C}_0)$,
\item $k=Gr(h_0)$,
\item $M=\mathrm{Im}(rb(R))$,
\item $\beta$ est la multi-dynamique graphique sur $\mathbf{G}$ d'ensemble de paramètres $M$ définie
pour tout sommet $S\in \dot{\mathbf{G}}$ par
\[S^\beta
=
\{(a_i)_{i\in I}\in S^{\alpha_0}\times \prod_{i\neq 0}(\Delta_i S)^{\alpha_i},
\forall i\neq 0, \tau_{i(\Delta_i S)}(a_i)=\delta_i(\tau_{0(S)}(a_0))\},\]
et pour toute arête $(e:S\rightarrow T)\in\overrightarrow{\mathbf{G}}$, tout paramètre $\mu\in M$ et tout état  $a=(a_i)_{i\in I}\in S^\beta$ par 
\[ 
e^\beta_\mu(a)=
\{
b\in T^\beta, 
\tau_{0(T)}(b_0)=e^{h_0}(\tau_{0(S)}(a_0))\,
\mathrm{et}\,
\exists (\mathfrak{a}_i)_{i\in I}\in  rb(R)^{-1}(\mu),
\forall i\in I, \mathfrak{a}_i\triangleright a_i, b_i
\}
\]
\item $\rho$ est la datation définie pour tout $a=(a_i)_{i\in I}\in S^\beta$ par $\rho_{(S)}(a)=\tau_{0(S)}(a_0)$. 
\end{itemize}

\subsection{Trois autres dynamiques ouvertes graphiques engendrées par $\mathcal{F}$}\label{subsec trois reduc sc pour F}

Appliquant à la famille dynamique graphique $\mathcal{G}=Gr(\mathcal{F})$ la définition donnée dans notre article \cite{Dugowson:20150807}, section §4.2, respectivement de
 la dynamique graphique $[\mathcal{G}]_\mathrm{f}$ fonctionnellement engendrée par $\mathcal{G}$,
de la dynamique graphique $[\mathcal{G}]_\mathrm{s}$ souplement engendrée par $\mathcal{G}$,
et de la mono-dynamique graphique $[\mathcal{G}]_\mathrm{m}$  engendrée par $\mathcal{G}$,
nous obtenons à présent à partir de notre famille dynamique sous-catégorique $\mathcal{F}$, respectivement trois dynamiques graphiques ouvertes\footnote{La troisième n'étant ouverte qu'au sens large, puisque son paramètre ne peut prendre qu'une unique valeur (voir \cite{Dugowson:20150807}, section §4.2.5).} :
 ${[Gr(\mathcal{F})]_\mathrm{f}}$,
 ${[Gr(\mathcal{F})]_\mathrm{s}}$,
et ${[Gr(\mathcal{F})]_\mathrm{m}}$.

\subsection{Théorème de stabilité}

\begin{thm}\label{thm stabilite sc}[Stabilité des dynamiques sous-catégoriques] Pour toute famille dynamique sous-catégorique $\mathcal{F}$, chacune des dynamiques graphiques ouvertes $[Gr(\mathcal{F})]_{\mathrm{p}}$, $[Gr(\mathcal{F})]_{\mathrm{f}}$, $[Gr(\mathcal{F})]_{\mathrm{s}}$ et $[Gr(\mathcal{F})]_{\mathrm{m}}$ est  sous-catégorique sur le moteur synchronisateur de $\mathcal{F}$.
\end{thm}

\begin{rmq}
Dans la preuve qui suit, et afin d'alléger l'écriture, on considère $\tau_i$ comme une application $st(\alpha_i)\rightarrow st(h_i)$ plutôt que comme une famille d'applications $(\tau_{i(S)}:S^{\alpha_i}\rightarrow S^{h_i})_{S\in \dot{\mathbf{C}_i}}$, de sorte que l'on écrira par exemple $\tau_{0}(a_0)$ plutôt que $\tau_{0(S)}(a_0)$. Même remarque pour les $\delta_i$. 
\end{rmq}

\paragraph{Preuve du théorème \ref{thm stabilite sc}.}
Reprenons les notations de la section \ref{subsec dyna graph primo-engendree par famille dysc} ci-dessus. Pour établir que 
la dynamique graphique ouverte
\[[Gr(\mathcal{F})]_{\mathrm{p}}=({\rho:(\beta:Gr(\mathbf{C}_0)\longrightarrow Gr(\mathbf{P}^{\underrightarrow{M}}))}\looparrowright Gr(h_0))\]
est sous-catégorique sur le moteur $\mathbf{C}_0$ de $\mathcal{F}$, l'horloge $h_0$ étant par définition catégorique sur $\mathbf{C}_0$ et la sous-catégoricité de $\rho$ étant, d'après le lemme \ref{lm dygrouv vers dysco}, une conséquence immédiate de la sous-catégoricité de $\beta$, il nous suffit de vérifier ce dernier point. Puisque $\beta$ est par construction une multi-dynamique graphique, $S^\beta_\mu=S^\beta$ ne dépend pas de $\mu\in M$. Reste donc à vérifier que, pour tout $\mu\in M$, les deux propriétés suivantes sont satisfaites :
\begin{enumerate}
\item  pour tout objet $S\in\dot{\mathbf{C}_0}$, 
\[(Id_S)^{\beta_\mu}\subset Id_{S^{\beta}},\]
\item  pour tout couple $(f,g)\in \overrightarrow{\mathbf{C}_0}^2$ de flèches composables
\[
(g\circ f)^{\beta_\mu} \subset g^{\beta_\mu}\odot f^{\beta_\mu}.
\]
\end{enumerate} 
Vérifions le premier point. Soit donc $\mu\in M$, $S\in\dot{\mathbf{C}_0}$, et $a=(a_i)_{i\in I}\in S^\beta$. Montrons que 
\[(Id_S)^{\beta_\mu}(a)\subset \{a\}.\]
Pour cela, supposant $(Id_S)^{\beta_\mu}(a)$ non vide, considérons un élément quelconque $b=(b_i)_{i\in I}\in (Id_S)^{\beta_\mu}((a_i)_{i\in I})$. Par définition de $\beta_\mu$, on a $b\in S^\beta$ qui vérifie
\[
\tau_{0(S)}(b_0)=(Id_S)^{h_0}(\tau_{0(S)}(a_0)).
\]
Mais puisque $h_0$ est une horloge catégorique, on en déduit $\tau_0(b_0)=\tau_0(a_0)$. Par ailleurs, toujours par définition de $\beta_\mu$, il existe  $(\mathfrak{a}_i)_{i\in I}\in  rb(R)^{-1}(\mu)$ telle que
\[
\forall i\in I, \mathfrak{a}_i\triangleright a_i, b_i
.
\]
On a alors, pour tout $i\in I$,
\[
\mathfrak{a}_i\triangleright  b_i 
\Rightarrow 
b_i
=\mathfrak{a}_i(\tau_i(b_i))
=\mathfrak{a}_i(\delta_i(\tau_0(b_0)))
=\mathfrak{a}_i(\delta_i(\tau_0(a_0)))
=\mathfrak{a}_i(\tau_i(a_i))
=a_i.
\]
Donc, on a bien soit $(Id_S)^{\beta_\mu}(a)=\emptyset$, soit $(Id_S)^{\beta_\mu}(a)=\{a\}$.

Vérifions finalement la seconde propriété.
Soit donc ${({(S\stackrel{f}{\rightarrow} T)},(T\stackrel{g}{\rightarrow} U))\in \overrightarrow{\mathbf{C}_0}^2}$ un couple de flèches composables, $\mu\in M$,  et $a=(a_i)_{i\in I}\in S^\beta$. Nous devons montrer que
\[
(g\circ f)^{\beta_\mu} (a)
\subset
(g^{\beta_\mu}\odot f^{\beta_\mu})(a).
\]
Supposons donc $(g\circ f)^{\beta_\mu} (a)$  non vide, et considérons un élément quelconque $c= (c_i)_{i\in I}\in (g\circ f)^{\beta_\mu} (a)$. Soit $(\mathfrak{a}_i)_{i\in I}\in  rb(R)^{-1}(\mu)$ une famille de parties externes de réalisations telle que
\[\forall i\in I, \mathfrak{a}_i\triangleright a_i, c_i,\] famille qui existe par définition de $(g\circ f)^{\beta_\mu}$.

Posons $t_0=\rho(a)=\tau_0(a_0)\in S^{h_0}$, $t_1=f^{h_0}(t_0)\in T^{h_0}$ et $t_2=g^{h_0}(t_1)\in U^{h_0}$. L'horloge $h_0$ étant catégorique, on a $t_2=(g\circ f)^{h_0}(t_0)$
et, par définition de $(g\circ f)^{\beta_\mu}(a)$, on a $(g\circ f)^{h_0}(t_0)=\tau_0(c_0)$, de sorte que $t_2=\tau_0(c_0)=\rho(c)$.

Posons ensuite $b=(b_i)_{i\in I}$ avec $b_0=\mathfrak{a}_0(t_1)\in T^{\alpha_0}$ et, pour tout $i\in I\setminus\{0\}$, 
\[ 
b_i=\mathfrak{a}_i(\delta_i(t_1)).
\]

Montrons que $b\in f^{\beta_\mu} (a)$. D'abord, on a bien $b\in T^\beta$ puisque, pour tout $i\neq 0$, on a premièrement
\[t_1\in T^{h_0}
\Rightarrow 
\delta_i(t_1)\in (\Delta_i T)^{h_i}
\Rightarrow 
\mathfrak{a}_i(\delta_i(t_1))\in (\Delta_i T)^{\alpha_i},\]
et deuxièmement
\[\tau_i(b_i)=\tau_i(\mathfrak{a}_i(\delta_i(t_1)))
=\delta_i(t_1)=\delta_i(\tau_0(b_0)) \]
puisqu'en effet $b_0=\mathfrak{a}_0(t_1) \Rightarrow t_1=\tau_0(b_0)$.
On en déduit au passage que $\rho(b)=t_1$. 
Ensuite, par construction même, on a $\forall i\in I, \mathfrak{a}_i\triangleright a_i, b_i$. On a donc bien $b\in f^{\beta_\mu} (a)$.

Montrons à présent que $c\in g^{\beta_\mu} (b)$. D'abord, on a bien $c\in U^\beta$, puisque $(g\circ f)^{\beta_\mu} (a)\subset U^\beta$. Ensuite, on a $t_2=\tau_0(c_0)=(g\circ f)^{h_0}(t_0)=g^{h_0}(t_1)$ d'où $\tau_0(c_0)=g^{h_0}(\tau_0(b_0))$. Enfin, la même famille de parties externes de réalisations $(\mathfrak{a}_i)_{i\in I}$ vérifie évidemment
\[
\forall i\in I, \mathfrak{a}_i\triangleright b_i, c_i.
\]
On a donc bien  $c\in g^{\beta_\mu} (b)$.

Finalement,  
\[
b\in f^{\beta_\mu} (a)
\,\mathrm{et}\,
c\in g^{\beta_\mu} (b)
\Rightarrow
c\in (g^{\beta_\mu} \odot f^{\beta_\mu})(a),\]
ce qu'il s'agissait d'établir.

Ayant ainsi établi que la dynamique primo-engendrée $[Gr(\mathcal{F})]_{\mathrm{p}}$ est sous-catégorique, le fait que les dynamiques graphiques $[Gr(\mathcal{F})]_{\mathrm{f}}$, $[Gr(\mathcal{F})]_{\mathrm{s}}$ et $[Gr(\mathcal{F})]_{\mathrm{m}}$ le soient également est une conséquence immédiate de la construction donnée aux sections \ref{subsubsec reduc param mdsc} et \ref{subsubsec reduc param mdsco}, elle-même fondée sur la proposition \ref{prop union dyna sous-cat}.
\begin{flushright}$\square$\end{flushright} 

\subsubsection{Dynamiques sous-catégoriques engendrées par $\mathcal{F}$}
Le théorème de stabilité des dynamiques sous-catégoriques permet de poser la définition suivante.

\begin{df}\label{df dynamiques sc engendrees} Soit $\mathcal{F}$ une famille dynamique sous-catégorique. On appelle \emph{dynamique ouverte primo-engendrée par $\mathcal{F}$}, et on note $[\mathcal{F}]_{\mathrm{p}}$  la dynamique sous-catégorique ouverte sur le moteur synchronisateur de $\mathcal{F}$ qui, en tant que dynamique graphique,  coïncide avec $[Gr(\mathcal{F})]_{\mathrm{p}}$. On définit de même, toujours sur le moteur synchronisateur de $\mathcal{F}$,  les dynamiques sous-catégoriques ouvertes $[\mathcal{F}]_{\mathrm{f}}$, $[\mathcal{F}]_{\mathrm{s}}$ et     $[\mathcal{F}]_{\mathrm{m}}$ :
\begin{itemize}
\item la \emph{dynamique ouverte $[\mathcal{F}]_{\mathrm{f}}$, fonctionnellement engendrée par $\mathcal{F}$}  coïncide avec $[Gr(\mathcal{F})]_{\mathrm{f}}$,
\item la \emph{dynamique ouverte $[\mathcal{F}]_{\mathrm{s}}$, souplement engendrée par $\mathcal{F}$}  coïncide avec $[Gr(\mathcal{F})]_{\mathrm{s}}$,
\item la \emph{mono-dynamique scandée $[\mathcal{F}]_{\mathrm{m}}$, mono-engendrée par $\mathcal{F}$} coïncide avec $[Gr(\mathcal{F})]_{\mathrm{f}}$.
\end{itemize}
\end{df}

\begin{rmq} Reprenant les constructions de la section §4.2 de \cite{Dugowson:20150807}, on voit qu'il est bien sûr possible de combiner les diverses relations d'équivalence sur l'ensemble des paramètres qui définissent les divers modes d'engendrements paramétriques --- ces relations d'équivalence étant elles-mêmes définies par le  choix de \og tas\fg\, au sein des paramètres de chacune des composantes de la famille dynamique considérée, tas dont il possible de prendre l'union, l'intersection, etc... --- pour produire de nouveaux modes d'engendrements paramétriques, l'idée étant de trouver le juste équilibre entre l'ouverture excessive de ${[\mathcal{F}]_\mathrm{p}}$, offrant des paramètres en réalité souvent inutilisables, et la fermeture complète sur elle-même de ${[\mathcal{F}]_\mathrm{m}}$, qui n'offre plus aucune prise à l'interaction avec d'autres dynamiques.
\end{rmq}

\subsection{Dynamiques catégoriques et stabilité}\label{subsec instab cat}


\begin{df} On appelle \emph{famille dynamique catégorique} une famille dynamique sous-catégorique dont toutes les composantes sont catégoriques. 
\end{df}

\begin{df} Une famille dynamique catégorique $\mathcal{F}$ est dite \emph{primo-stable} (respectivement : \emph{fonctionnellement stable}, \emph{souplement stable}, \emph{mono-stable}) si la dynamique sous-catégorique ouverte $[\mathcal{F}]_{\mathrm{p}}$ (respectivement : $[\mathcal{F}]_{\mathrm{f}}$, $[\mathcal{F}]_{\mathrm{s}}$, $[\mathcal{F}]_{\mathrm{m}}$) est catégorique. 
\end{df}

\section{Conclusion}\label{section conclusion DySCO}

Tout comme pour l'article précédent \cite{Dugowson:20150807}, notre objectif n'était pas ici de donner une intuition des notions présentées  mais d'écrire les définitions qui, ultérieurement, permettront de présenter les exemples qui  montreront, nous l'espérons,  tout leur intérêt. Deux phénomènes en particulier ne manqueront pas alors d'attirer l'attention.

D'une part, 
et c'est là un phénomène très intéressant aussi d'un point de vue philosophique, dans certains cas, l'interaction de dynamiques déterministes peut produire une dynamique non déterministe, du moins dans le cas de l'engendrement fonctionnel.  Décrivons en quelques mots le principe d'un tel phénomène  : considérant deux systèmes en interaction, et donnant  à chacun d'eux comme contrainte déterminante --- \og déterminante\fg\, en tout cas selon le point de vue de l'engendrement fonctionnel --- de se comporter exactement comme l'autre, alors la dynamique ainsi engendrée est très largement indéterministe. 

D'autre part, comme nous l'avons annoncé d'emblée, il y a le fait --- qui est à l'origine de notre intérêt pour les dynamiques sous-catégoriques --- que l'interaction de dynamiques catégoriques n'engendre pas toujours une dynamique catégorique, ce dont il est très facile de donner des exemples\footnote{\label{footnote principe exemple cat instable} Nous avons donné le principe d'un tel exemple dans notre conférence \cite{Dugowson:20150506} du 6 mai 2015.}. En fait, au moins lorsqu'on s'en tient à l'usage des horloges existentielles\footnote{Voir \cite{Dugowson:201112} et \cite{Dugowson:201203}, section §3.4.3.},
cette instabilité des dynamiques catégoriques ne s'observe pas pour les moteurs catégoriques classiques ($\mathbf{N}$, $\mathbf{Z}$, $\mathbf{R}_+$ et $\mathbf{R}$): dans ces cas là, en effet, les horloges existentielles sont totalement ordonnées, et il est toujours possible, en particulier, de recoller une solution définie sur un intervalle $[a,b]$ avec une solution définie sur un intervalle $[b,c]$. Pour certains moteurs catégoriques, les familles dynamiques catégoriques à horloges existentielles seront ainsi toujours stables et, avec ces horloges, c'est donc uniquement lorsqu'entrent en jeu des temporalités non ordinaires --- par exemple cycliques, arborescentes, multi-dimensionnelles, etc... --- que l'instabilité des dynamiques catégorique est susceptible de se manifester. Les conditions de stabilité des interactions entre dynamiques catégoriques constituent ainsi l'une des explorations qu'il s'agira de poursuivre. 

Enfin, l'étape suivante naturelle pour toutes ces recherches concerne les dynamiques catégoriques (ou sous-catégoriques) ouvertes \emph{connectives}, et leurs interactions.

\bibliographystyle{plain}




\begin{thebibliography}{}

\end{thebibliography}


\begin{thebibliography}{1}

\bibitem{Dugowson:201012}
Stéphane Dugowson.
\newblock On connectivity spaces.
\newblock {\em Cahiers de {T}opologie et {G}éométrie {D}ifférentielle
  {C}atégoriques}, LI(4):282--315, 2010.
\newblock \url{http://hal.archives-ouvertes.fr/hal-00446998/fr}.

\bibitem{Dugowson:201112}
Stéphane Dugowson.
\newblock Introduction aux dynamiques catégoriques connectives, décembre
  2011.
\newblock \url{{http://hal.archives-ouvertes.fr/hal-00654494/fr/}}.

\bibitem{Dugowson:201203}
Stéphane Dugowson.
\newblock {\em Dynamiques connectives (Une introduction aux notions connectives
  : espaces, représentations, feuilletages et dynamiques catégoriques)}.
\newblock \'{E}ditions Universitaires Européennes, 2012.

\bibitem{Dugowson:20130521}
Stéphane Dugowson.
\newblock Influences dynamiques catégoriques, mai 2013.
\newblock Exposé lors du colloque : \emph{Topologie, effectivité,
  interactivité} organisé par Supméca et Innovaxiom (21 mai 2013). Vidéo accessible en ligne à
  l'adresse : {\url{{https://sites.google.com/site/dugowsonrecherche/jms/jms-4}}}.


\bibitem{Dugowson:201505}
Stéphane Dugowson.
\newblock Structure connective des relations multiples, mai 2015.
\newblock \url{https://hal.archives-ouvertes.fr/hal-01150262}.

\bibitem{Dugowson:20150506}
Stéphane Dugowson.
\newblock {Possibilités de principe et possibilités effectives dans les
  dynamiques catégoriques}, mai 2015.
\newblock Exposé lors du colloque : \emph{Théorie des catégories, dynamiques
  anciennes et nouvelles, mathématiques et philosophie}. Paris Diderot (5 et 6
  mai 2015). Vidéos accessibles en ligne sur le site CLE, à l'adresse :
  {\url{https://sites.google.com/site/logiquecategorique/autres-seminaires/DYAN/20150506Dugowson}}.
  
  
  
\bibitem{Dugowson:20150807}
Stéphane Dugowson.
\newblock {Interaction des dynamiques graphiques ouvertes (définitions)},
  {août} 2015.
\newblock \url{https ://hal.archives-ouvertes.fr/hal-01177450}.



\end{thebibliography}

\tableofcontents

\end{document}